\documentclass[3pt]{article}

\usepackage{caption}
\usepackage{graphicx}\usepackage{epstopdf}\usepackage{color}
\usepackage{amsmath}\usepackage{amsfonts}\usepackage{amssymb}
 \usepackage{amsthm}
\newtheorem{theorem}{Theorem}[section]
\newtheorem{lemma}[theorem]{Lemma}

\numberwithin{equation}{section}

\begin{document}

\title{An $H^m$-conforming spectral element method on 
multi-dimensional domain and its application  to  transmission eigenvalues}

\author{ {Jiayu Han, Yidu Yang\thanks{Corresponding Author}} \\\\%
{\small School of Mathematics and Computer Science, }\\{\small
Guizhou Normal University,  Guiyang,  $550001$,  China}\\{\small
hanjiayu126@126.com, ydyang@gznu.edu.cn}}
\date{~}
\pagestyle{plain} \textwidth 145mm \textheight 215mm \topmargin 0pt
\maketitle

\indent{\bf\small Abstract~:~} {\small 
In this paper we develop an $H^m$-conforming ($m\ge1$) spectral element method on multi-dimensional domain associated with the partition into multi-dimensional rectangles. We construct a set of basis functions on the interval $[-1,1]$ that is  made up of  the generalized Jacobi polynomials (GJPs) 
 and the nodal basis functions.
So the basis functions on multi-dimensional rectangles consist of the tensorial product of the basis functions on the interval $[-1,1]$.
Then we construct the spectral element interpolation operator and   prove the associated interpolation error estimates. Finally we apply the $H^2$-conforming spectral element method  to the Helmholtz transmission
eigenvalues that is a hot topic in the field  engineering and  mathematics.
%
}\\
\indent{\bf\small Keywords~:~\scriptsize} {\small Spectral element method,
Multi-dimensional domain, Interpolation error estimates, Transmission eigenvalues}
\section{Introduction}
\label{intro}

Spectral  method is an efficient method in scientific and engineering  computations which can provide superior accuracy for the solution of partial differential equations in fluid dynamics \cite{shen,canuto1}.
But spectral method lacks the domain flexibility. So the spectral element method is developed to overcome  this defect.  Up to now the spectral element method
have attracted more and more scholars's attention. Guo and Jia \cite{guo} studied the quadrilateral spectral method and extended it to $H^1$-conforming spectral element method for polygons. Shen et al. \cite{shen3}  provided an $H^1$-conforming  spectral element method by constructing directly the modal basis functions on the triangle while Samson  et al. \cite{samson} builds a new $H^1$-conforming spectral element method  using the basis on the triangle 
by the rectangle-triangle mapping. Yu and Guo \cite{yu} developed an $H^2$-conforming spectral element method with  rectangular partition  in two dimension.

The orthogonal Jacobi polynomials $\widehat{J}_i^{\alpha,\beta}(\widehat{x})$ ($\widehat{x}\in I:=[-1,1]$, $\alpha,\beta>-1$) weighted with $w^{\alpha,\beta}(\widehat{x})$ $:=(1-\widehat{x})^\alpha(1+\widehat{x})^\beta$ are usually adopted to form the modal basis functions in spectral method and spectral element method.  Shen et al. \cite{shen} extend these polynomials to the generalized case, namely the  GJPs  for $\alpha,\beta\in\mathbb R$. It is worth indicating  an important property of GJPs that
they together with their  first few  order derivatives vanish  at the endpoints $\pm1$. So Shen et al.
use them as a set of basis functions in $H^m_0(I)$ $(m\ge1)$ as well as apply them to the general order PDEs. Using the GJPs one can easily construct the tensional basis functions in $H^m_0$ for spectral method on multi-dimensional rectangle. Canuto et al. mentioned in section 8.5 of the book \cite{canuto1} a type of modal boundary-adapted basis functions on [-1,1].    They consist of the modal basis, viewed as a compact combination of Legendre polynomials,  and the nodal basis (without derivatives) at $\pm1$ so that    one can easily establish $H^1$-conforming spectral element approximation. Note that these modal bases are no other than  the GJPs $\widehat{J}_i^{-1,-1}(\widehat{x})$. Using the similar way Yu and Guo \cite{yu} developed an $H^2$-conforming spectral element method with  rectangular partition   coupled with an error analysis.
This motivates us to extend this situation to $H^m$-conforming spectral elements on multi-dimensional domain. 

In this paper, we aim to develop an $H^m$-conforming spectral element method on multi-dimensional domain which is the same as the one in \cite{yu} for the case $m=2$ in two dimension.
    We construct a set of basis functions on the interval $[-1,1]$ that is  made up of  GJPs, which is regarded as the bubble functions,  and the nodal basis functions.
So the basis functions on multi-dimensional rectangular element consist of the tensorial product of the basis functions on the interval $[-1,1]$ by an affine mapping.
Then we construct the spectral element interpolation operator and   prove the associated interpolation error estimates.
 Finally we shall apply the $H^2$-conforming spectral element method presented in this paper to the Helmholtz transmission eigenvalue problem that is a quadratic eigenvalue problem arising in  inverse scattering theory for an inhomogeneous medium \cite{cakoni2,colton1}. 
   In recent years, the numerical methods of the transmission eigenvalue problem are   hot topics in the field of engineering and computational mathematics (see \cite{colton2,an2,ji,sun1,ji2,cakoni4,yang1}). Among them   \cite{an2} studied the spectral methods on the rectangle.
   But to our knowledge the above works do not involve   spectral element method with $d$-dimensional rectangular partition  ($d=2,3$). In this paper   we adopt the $H^2$-conforming  method builded in \cite{yang1} to construct a spectral element approximation for  transmission eigenvalues. Our theoretical analysis and numerical results show that the $H^2$-conforming spectral element method can obtain the
   transmission eigenvalues of high accuracy  numerically.

   \section{An $H^m$-conforming spectral element method}
In this section, we shall discuss an $H^m$-conforming spectral element method on $d$-dimensional domain $D$ ($d\ge1$).
We associate $D$ with   a  sequence of rectangular partitions $\{\pi_h\}_{h>0}$ into elements $\kappa$.
First of all, we   consider the construction of the
bases on one-dimensional standard interval $I=[-1,1]$  containing nodal bases and   modal bases. Before presentation, we use the notation $P_N(K)$ to denote the polynomial space of degree $\le N$ in each variable on $K$.

 One  need to construct $2m$ nodal basis functions $\widehat\phi_{j}(\widehat x)$ ($j=0,\cdots,2m-1$) for the polynomial space $P_{2m-1}(I)$ satisfying
$$\partial_{\widehat x}^j\widehat\phi_{i}(-1)=\partial_{\widehat x}^j\widehat\phi_{i+m}(1)=\delta_{i,j}, (i,j=0,\cdots,m-1).$$

Then we shall construct the modal bases  on $I$ which are actually the polynomial bubble functions on $I$.
Let $\widehat{J}_j^{\alpha,\beta}(\widehat{x})$ be the Jacobi polynomials which are orthogonal with respect to the  weight function $\widehat{\omega}^{\alpha,\beta}(\widehat{x})=(1-\widehat{x})^\alpha(1+\widehat{x})^\beta$ ($\alpha,\beta>-1$) on $I$:

\begin{eqnarray*}
 \int_{-1}^1 \widehat{J}_i^{\alpha,\beta}(\widehat{x})\widehat{J}_j^{\alpha,\beta}(\widehat{x})\omega^{\alpha,\beta}(\widehat{x})=\gamma_j^{\alpha,\beta}\delta_{i,j}.
\end{eqnarray*}
where $\gamma_j^{\alpha,\beta}=\frac{2^{\alpha+\beta+1}\Gamma(j+\alpha+1)\Gamma(j+\beta+1)}{(2j+\alpha+\beta+1)j!\Gamma(j+\alpha+\beta+1)}$.
 The  GJPs  are defined by
 \begin{eqnarray*}
\widehat{J}_j^{\alpha,\beta}(\widehat{x}) =\left\{\begin{tabular}{ll}
                                      $(1-\widehat{x})^{-\alpha}(1+\widehat{x})^{-\beta}\widehat{J}_{j-j_0}^{-\alpha,-\beta}(\widehat{x})$,&$\alpha,\beta\le-1$,\\
                                      $(1-\widehat{x})^{-\alpha}\widehat{J}_{j-j_0}^{-\alpha,\beta}(\widehat{x})$,&$\alpha\le-1,\beta>-1$,\\
                                      $(1+\widehat{x})^{-\beta}\widehat{J}_{j-j_0}^{\alpha,-\beta}(\widehat{x})$,&$\alpha>-1,\beta\le-1$,\\
                                   \end{tabular}\right.
\end{eqnarray*}
where $j\ge j_0$ with $j_0=-(\alpha+\beta),-\alpha,-\beta$ for the above three cases, respectively.

 Here we fix $\alpha=\beta=-m$ then the GJPs   $\{\widehat{J}^{-m,-m}_{j}(\widehat x)\}_{j\ge 2m}$ 
 satisfy
$$\int_{-1}^1 \widehat{J}^{-m,-m}_{i}(\widehat x)\widehat{J}^{-m,-m}_{j}(\widehat x)w^{-m,-m}(\widehat x)d\widehat x=\gamma^{m,m}_{j-2m} \delta_{i,j}.$$
An attractive property of GJPs is that
\begin{eqnarray*}
\partial_{\widehat x}^j {\widehat{J}}_i^{-m,-m}(\pm1)=0,  ~j=0,1,\cdots,m-1,~i\ge 2m.
\end{eqnarray*}
In addition, GJPs can be represented as a compact combination of Legendre polynomials (see Lemma 6.1 in \cite{shen} and Remark 1) which is convenient for computations. So we adopt them to set the bubble functions on $I$
$$\widehat\phi_{j}(\widehat x)=\widehat{J}^{-m,-m}_{j}(\widehat x),~j=2m,2m+1,\cdots,N.$$
 It is known that $\big\{\widehat\phi_{j}\big  \}_{j=2m}^{N}$ is a set of basis functions of $P_N^0(I)\subset H_0^m(I)$ (see \cite{shen}).
Hence $\big\{\widehat\phi_{j}\big  \}_{j=0}^{N}$ constitutes a set of basis functions of $P_N(I)$.
Next we consider the case of an arbitrary interval $[a,b]$.
We define $$\phi_{j}(x) =(\frac{b-a}{2})^j\widehat\phi_{j}(\widehat x)~and~\phi_{j+m}(x) =(\frac{b-a}{2})^j\widehat\phi_{j+m}(\widehat x)~for~0\le j\le m-1$$
$$\phi_{j}(x) =\widehat\phi_{j}(\widehat x)~for~2m\le j\le N$$
in terms with the linear transformation $$x=\frac{b-a}{2}\widehat x+\frac{b+a}{2}$$ then $\big\{ \phi_{j}\big  \}_{j=0}^{N}$ constitutes a set of basis functions of $P_N([a,b])$.

Now we consider the bases for the arbitrary element $\kappa:=[a_1,b_1]\times\cdots\times[a_d,b_d]\subset R^d$.
A natural choice of the bases on $\kappa$ is the tensor product of one-dimensional basis functions.
We define the linear transformation: $x_i=\frac{b_i-a_i}{2}\widehat x_i+\frac{b_i+a_i}{2}$ ($i=1,\cdots,d$).
Based on the previous discussion for one dimension, one can use $\big\{\prod\limits_{i=1}^d\phi_{j_i}(x_i) \big  \}_{j_1,\cdots,j_d=0}^{N}\subset  {P}_N(\kappa)$ as a set of basis functions for the element $\kappa$.
For reading conveniently, we   classify these basis functions on $\kappa$ as follows:

Nodal basis functions: $\prod\limits_{i=1}^d\phi_{j_i}(x_i),  j_1,\cdots,j_d=0,\cdots,2m-1.$

$q$-face basis functions($1\le q\le d-1$): For   any rearranged sequence $\{i_{l}\}_{l=1}^d$ of $\{i\}_{i=1}^d$,  define $ \prod\limits_{i=1}^d\phi_{j_i}(x_i)$,    $j_{i_1},\cdots,~j_{i_{d-q}}=0,\cdots,2m-1,~j_{i_{d-q+1}},\cdots,j_{i_{d}}\ge 2m$.

Element bubble basis functions: $\prod\limits_{i=1}^d\phi_{j_i}(x_i),  j_1,\cdots,j_d\ge 2m.$

One can easily verify the $H^m$-conformity for basis functions between the adjacent element $\kappa_1$ and $\kappa_2$.  We consider only the case that $\kappa_1$ and $\kappa_2$ share the common $(d-1)$-face $\partial \kappa_1\cap\partial\kappa_2=a\times[a_2,b_2]\times\cdots\times[a_d,b_d]$.

For $s=0,\cdots,m-1$, we have
\begin{eqnarray*}
&&\partial_{x_1}^{s}\phi_{j_1}|_{\kappa_1}(a)=\partial_{x_1}^{s}\phi_{j_1}|_{\kappa_2}(a),~ j_{1} =0,\cdots,2m-1,\\
&&\partial_{x_1}^{s}\phi_{j_1}|_{\kappa_1}(a)=\partial_{x_1}^{s}\phi_{j_1}|_{\kappa_2}(a)=0,~ j_{1} \ge 2m.
\end{eqnarray*}

For $s=0,\cdots,m-1$ and $i=2,\cdots,d$, we have with $x_i\in [a_i,b_i]$
\begin{eqnarray*}
&& \partial_{x_i}^{s}\phi_{j_i}|_{\kappa_1}(x_i)=\partial_{x_i}^{s}\phi_{j_i}|_{\kappa_2}(x_i)=(\frac{2}{b_i-a_i})^{s-j_i}\partial_{\widehat x_i}^{s}\widehat{\phi}_{j_i}(\widehat x_i) ,~ j_{i}= 0,\cdots,m-1,\\
&& \partial_{x_i}^{s}\phi_{j_i}|_{\kappa_1}(x_i)=\partial_{x_i}^{s}\phi_{j_i}|_{\kappa_2}(x_i)=(\frac{2}{b_i-a_i})^{s-j_i+m}\partial_{\widehat x_i}^{s}\widehat{\phi}_{j_i}(\widehat x_i) ,~ m\le j_{i}\le 2m-1,\\
&& \partial_{x_i}^{s}\phi_{j_i}|_{\kappa_1}(x_i)=\partial_{x_i}^{s}\phi_{j_i}|_{\kappa_2}(x_i)=(\frac{2}{b_i-a_i})^s\partial_{\widehat x_i}^{s}\widehat{J}^{-m,-m}_{j_i}(\widehat x_i) ,~ j_{i} \ge 2m .
\end{eqnarray*}

Therefore, the   bases $\prod\limits_{i=1}^d\phi_{j_i}(x_i),  (0\le j_1,\cdots,j_d \le N)$ together with their derivatives of order $\le m-1$ are equal on $\partial \kappa_1\cap\partial\kappa_2$.

In what follows, we mainly introduce some interpolation operators that will be used in the argument afterwards.

Introduce the  interpolation operator  $\widehat\Pi_i^1$ ($1\le i\le d$): $$(\widehat{\Pi}_i^1 \widehat v)(\widehat x_i)=\sum_{j=0}^{m-1}((\partial^j_{\widehat x_i}\widehat v)(-1)\widehat\phi_j(\widehat x_i)+(\partial^j_{\widehat x_i} \widehat v)(1)\widehat\phi_{j+m}(\widehat x_i)),$$
and the orthogonal projector $\widehat\Pi_i^2$ from $H^m_0(I)$ to $P_N^0(I)=P_N(I)\cap H^m_0(I)$: $$\int_{-1}^1\partial_{\widehat x_i}^m(\widehat{\Pi}_i^2 \widehat v(\widehat x_i)-\widehat v(\widehat x_i))\partial_{\widehat x_i}^m\widehat{v}_N(\widehat x_i)=0,~\forall \widehat{v}_N\in P_N^0(I).$$
Define $ (\Pi_i^1 v)(  x_i)=(\widehat\Pi_i^1\widehat v)(\widehat x_i)$
and
 $ (\Pi_i^2 v)(  x_i)=(\widehat\Pi_i^2\widehat v)(\widehat x_i)$ with $v(  x_i)=\widehat v(\widehat x_i)$. Then it is obvious that $\Pi_i^2$
is an orthogonal projector   from $H^m_0([a_i,b_i])$ to $P^0_N([a_i,b_i])$
and $$( {\Pi}_i^1   v)(  x_i)=\sum_{j=0}^{m-1}((\partial^j_{  x_i}   v)(-1) \phi_j(  x_i)+(\partial^j_{  x_i}   v)(1) \phi_{j+m}(  x_i)).$$
Define one-dimensional interpolation operator $\widehat{I}_{i}:C^{m-1}(I)\rightarrow P_N(I)$ and $(\widehat I_{i}\widehat v)( {\widehat x_i})=\widehat v( \widehat{x}_i)$ ($i=1,\cdots,d$)
satisfying
$$(\widehat I_{i}\widehat v)(\widehat x_i)= (\widehat\Pi_i^1\widehat v+\widehat \Pi_i^2\circ(\mathbf{I}-\widehat\Pi_i^1) \widehat v)(\widehat x_i),$$
where $\mathbf{I}$ is the identity operator.

Define the function $v(x_i)=\widehat v (\widehat x_i)$ and $ {I}_{i}:C^{m-1}([a_i,b_i])\rightarrow P_N([a_i,b_i])$ satisfying $(  I_{i}  v)( {  x_i})=  v(  {x}_i),\forall x_i\in [a_i,b_i]$ and
$$(  I_{i}  v)(  x_i)=( \Pi_i^1 v +  \Pi_i^2\circ(\mathbf{I}-\Pi_i^1) v)(  x_i).$$
It is obvious that $(  I_{i}  v)(  x_i)=(\widehat I_{i}\widehat v)(\widehat x_i)$.

\indent Let $H^{s}(K)$ be a Sobolev space with norm $\|\cdot\|_{s,K}$ for a given $K  \subseteq D$ and
we shall omit the subscript $K$ if $K=D$. Hereafter in this paper, we use the symbols  $x \lesssim y$ to mean $x \le Cy$ for a constant $C$ that is independent of the mesh size and the degree of piecewise polynomial space and may be different at different occurrences. Now we start with the one-dimensional interpolation error estimates.

\begin{lemma}
Assume $\widehat v\in   H^t(I)$ ($m\le t\le N+1$) then there holds for $0\le s\le m$
\begin{eqnarray*}
\|\widehat I_{i}\widehat v-\widehat v\|_{s,I}
  \lesssim (1/N)^{t-s}\|{\widehat v}\|_{t,I}.
\end{eqnarray*}
\begin{proof} We consider only the case of the integers $s$ and $t$ since the left case can be derived by the operator interpolation theory.
From   Theorem 6.1 in \cite{shen} we know if $\widehat v\in H^m_0(I)\cap H^t(I)$ ($t\ge m$) then there holds for $0\le s\le m$
\begin{eqnarray*}
\|\widehat \Pi_i^2 \widehat v-\widehat v\|_{s,I}\lesssim (1/N)^{t-s}\|\widehat v\|_{t,I}.
\end{eqnarray*}
Note that $\widehat v-\widehat\Pi_i^1 \widehat v\in H^m_0(I)$ and   $\widehat\Pi_i^1(\widehat v-\widehat\Pi_i^1 \widehat v)=0$. Then
\begin{eqnarray*}
\|\widehat I_{i}\widehat v-\widehat v\|_{s,I}&=& \|\widehat I_{i}(\widehat v-\widehat\Pi_i^1 \widehat v)-(\widehat v-\widehat\Pi_i^1 \widehat v)\|_{s,I}\\
&=&\|\widehat \Pi_i^2\circ(\mathbf{I}-\widehat{\Pi}_i^1)(\widehat v-\widehat\Pi_i^1 \widehat v)-(\widehat v-\widehat\Pi_i^1 \widehat v)\|_{s,I}\\
&\le&\|\widehat \Pi_i^2\circ(\mathbf{I}-\widehat{\Pi}_i^1)(\widehat v-\widehat\Pi_i^1 \widehat v)-  (\mathbf{I}-\widehat{\Pi}_i^1)(\widehat v-\widehat\Pi_i^1 \widehat v) \|_{s,I}\\
 & \lesssim& (1/N)^{t-s}\|(\mathbf{I}-\widehat{\Pi}_i^1)(\widehat v-\widehat\Pi_i^1 \widehat v)\|_{t,I}\\
 & \lesssim& (1/N)^{t-s}\|\widehat v\|_{t,I}.
\end{eqnarray*}
This concludes the proof.
\end{proof}
\end{lemma}
Define the $d$-dimensional interpolation operator $\widehat{\mathbf{I}}_{N}:C^{m-1}(I^d)\rightarrow P_N(I^d)$ as $$\widehat{\mathbf{I}}_{N}=\widehat{ I}_{1}\circ \cdots\circ\widehat{ I}_{d}.$$
One can easily verify that  $(\widehat{\mathbf{I}}_{N} \widehat v)(\mathbf{\widehat x})=\widehat v(\mathbf{\widehat x})$ holds for any $\widehat v(\mathbf{\widehat x})\in P_N(I^d)$.

\begin{lemma} For any $\widehat v\in H^t(I^d)$ with $md\le t\le N+1$
\begin{eqnarray*}
\|\widehat{\mathbf{I}}_{N} \widehat v-\widehat v\|_{s,I^d}
 \lesssim  (1/N)^{t-s}\|\widehat v\|_{t,I^d},   ~ 0\le s\le m .
\end{eqnarray*}
\begin{proof} Similar to Lemma 2.1 we only consider the case when both $s$ and $t$ are integers. Let $d$  nonnegative integers $\alpha_1,\cdots,\alpha_d$ satisfy
$\sum_{i=1}^d\alpha_i=s$.
We obtain from Lemma 2.1 that
\begin{eqnarray*}
\|\partial_{x_1}^{\alpha_1}\cdots\partial_{x_d}^{\alpha_d}(\widehat{ I}_{1} \widehat v-\widehat v)\|_{0,I^d}
  \lesssim  (1/N)^{t-s}\|  \widehat v\|_{t,I^d} ,
\end{eqnarray*}
and
\begin{eqnarray*}
&&\|\partial_{x_1}^{\alpha_1}\partial_{x_2}^{\alpha_2}\cdots\partial_{x_d}^{\alpha_d}(\widehat{ I}_{1}-\mathbf{I})(\widehat{ I}_{2}\circ \cdots\circ\widehat{ I}_{d}\widehat v -\widehat v)\|_{s,I^d}\\
&&~~~\lesssim
 (1/N)^{m-\alpha_1}\| \partial_{x_1}^{m}\partial_{x_2}^{\alpha_2}\cdots\partial_{x_d}^{\alpha_d} (\widehat{ I}_{2}\circ \cdots\circ\widehat{ I}_{d} \widehat v  - \widehat v)\|_{0,I^d},
\end{eqnarray*}
where $\mathbf{I}$ is the identity operator.
Hence by the triangular inequality
\begin{eqnarray*}
&&\|\partial_{x_1}^{\alpha_1}\cdots\partial_{x_d}^{\alpha_d}(\widehat{\mathbf{I}}_{N} \widehat v-\widehat v)\|_{0,I^d}\\
&&~~~
\lesssim\|\partial_{x_1}^{\alpha_1}\cdots\partial_{x_d}^{\alpha_d}(\widehat{ I}_{1} \widehat v-\widehat v)\|_{0,I^d}+\|\partial_{x_1}^{\alpha_1}\cdots\partial_{x_d}^{\alpha_d}(\widehat{ I}_{2}\circ \cdots\circ\widehat{ I}_{d} \widehat v -\widehat v)\|_{0,I^d}\\
&&~~~~~~+\|\partial_{x_1}^{\alpha_1}\cdots\partial_{x_d}^{\alpha_d}(\widehat{ I}_{1}-\mathbf{I})(\widehat{ I}_{2}\circ \cdots\circ\widehat{ I}_{d}\widehat v -\widehat v)\|_{0,I^d}\\
&&~~~\lesssim (1/N)^{t-s}\| \widehat v\|_{t,I^d}+\|\partial_{x_1}^{\alpha_1}\cdots\partial_{x_d}^{\alpha_d}(\widehat{ I}_{2}\circ \cdots\circ\widehat{ I}_{d} \widehat v -\widehat v)\|_{0,I^d}\\
&&~~~~~~+(1/N)^{m-\alpha_1}\| \partial_{x_1}^{m}\partial_{x_2}^{\alpha_2}\cdots\partial_{x_d}^{\alpha_d} (\widehat{ I}_{2}\circ \cdots\circ\widehat{ I}_{d} \widehat v  - \widehat v)\|_{0,I^d}\\
&&~~~\lesssim (1/N)^{t-s}\| \widehat v\|_{t,I^d}+\|\partial_{x_1}^{\alpha_1} \cdots\partial_{x_d}^{\alpha_d}(\widehat{ I}_{3}\circ \cdots\circ\widehat{ I}_{d} \widehat v -\widehat v)\|_{0,I^d}\\
&&~~~~~~+(1/N)^{m-\alpha_2}\| \partial_{x_1}^{\alpha_1}\partial_{2}^{m}\cdots\partial_{x_d}^{\alpha_d} (\widehat{ I}_{3}\circ \cdots\circ\widehat{ I}_{d} \widehat v  - \widehat v)\|_{0,I^d}\\
&&~~~~~~+(1/N)^{m-\alpha_1}\| \partial_{x_1}^{m}\partial_{x_2}^{\alpha_2}\cdots\partial_{x_d}^{\alpha_d} (\widehat{ I}_{3}\circ \cdots\circ\widehat{ I}_{d} \widehat v  - \widehat v)\|_{0,I^d}\\
&&~~~~~~+(1/N)^{(m-\alpha_1)+(m-\alpha_2)}\| \partial_{x_1}^{m}\partial_{x_2}^{m}\cdots\partial_{x_d}^{\alpha_d} (\widehat{ I}_{3}\circ \cdots\circ\widehat{ I}_{d} \widehat v  - \widehat v)\|_{0,I^d}.
\end{eqnarray*}
Repeating the above argument we get
\begin{eqnarray*}
&&\|\partial_{x_1}^{\alpha_1}\cdots\partial_{x_d}^{\alpha_d}(\widehat{\mathbf{I}}_{N} \widehat v-\widehat v)\|_{0,I^d} \lesssim (1/N)^{t-s}\| \widehat v\|_{t,I^d}+\|\partial_{x_1}^{\alpha_1} \cdots\partial_{x_d}^{\alpha_d}(\widehat{ I}_{d} \widehat v -\widehat v)\|_{0,I^d}\\
&&~~~~~~+\sum_{k=1}^{d-1}\sum_{i_1<\cdots<i_k\le d-1}(\frac{1}{N})^{(m-\alpha_{i_1})+\cdots+(m-\alpha_{i_k})}\| \partial_{x_{i_1}}^{m}\cdots\partial_{x_{i_k}}^{m} \partial_{x_{i_{k+1}}}^{\alpha_{i_{k+1}}}\cdots\partial_{x_{i_{d-1}}}^{\alpha_{i_{d-1}}}\partial_{x_{d}}^{\alpha_{d}} \\
&&~~~~~~~~~(   \widehat{ I}_{d} \widehat v  - \widehat v)\|_{0,I^d}\\
&&~~~\lesssim(1/N)^{t-s}\| \widehat v\|_{t,I^d}.
\end{eqnarray*}
This ends this proof.
\end{proof}
\end{lemma}

 Likewise $ { {\mathbf{I}}}_{N}^\kappa:C^1(\kappa)\rightarrow P_N(\kappa)$ is defined as $$ {\mathbf{I}}_{N}^\kappa= { I}_{1}\circ \cdots\circ { I}_{d}.$$

In the end, we introduce the spectral element space
 \begin{eqnarray*}
S^{N,h}=\{v:v|_\kappa\in P_N(\kappa),~\forall\kappa\in\pi_h ~\mathrm{and}~ \partial_{x_i}^s  v~(0\le i\le d,0\le s\le m-1)  \\
 \mathrm{~are ~ continuous ~accross ~\partial \kappa_1\cap\partial \kappa_2~for~\kappa_1,\kappa_2\in\pi_h ~\mathrm{and}~ \partial\kappa_1\cap\partial \kappa_2\neq \Theta}\}.
 \end{eqnarray*}

  We   define the spectral element interpolation operator $ { {\mathbf{I}}}_{N,h}:C^1(D)\rightarrow S^{N,h}$ as
${ {\mathbf{I}}}_{N,h}|_\kappa={ {\mathbf{I}}}_{N}^\kappa$ for any $\kappa\in \pi_h.$
One can easily verify that  $( {\mathbf{I}}_{N,h}   v)(\mathbf{  x})=  v(\mathbf{  x})$ holds for any $  v(\mathbf{  x})\in P_N(\kappa)$ and
$(\mathbf{I}_{N,h}v)(\mathbf{x})=(\widehat{\mathbf{I}}_{N}\widehat v)(\mathbf{\widehat{x}})$ for any $\mathbf{x}\in\kappa$.


Using the scaling argument, we can easily derive the interpolation error estimate on the element $\kappa$ and the entire domain $D$.

\begin{lemma} For any $  v\in H^t(\kappa)$ with $md\le t\le N+1$
\begin{eqnarray*}
\| {\mathbf{I}}_{N,h}   v-  v\|_{s,\kappa}
\lesssim (h/N)^{t-s}\|  v\|_{t,\kappa},   ~ 0\le s\le m .
\end{eqnarray*}
\end{lemma}

\begin{theorem} For any $  v\in H^t(D)$ with $md\le t\le N+1$
\begin{eqnarray*}
\| {\mathbf{I}}_{N,h}   v-  v\|_{s,D}
\lesssim (h/N)^{t-s}\|  v\|_{t,D},   ~ 0\le s\le m .
\end{eqnarray*}
\end{theorem}

\indent{\bf Remark 1.}~ We consider the special case $m=2$.
The   nodal basis functions with respective to the endpoint -1 are respectively:
$$\widehat{\phi}_0(\widehat x)=\frac{(\widehat x-1)^2(\widehat x+2)}{4}, ~\widehat\phi_1(\widehat x)=\frac{(\widehat x-1)^2(\widehat x+1)}{4}.$$
Meanwhile,
the    nodal basis functions with respective to the endpoint 1 are respectively:
$$\widehat\phi_2(\widehat x)=-\frac{(\widehat x+1)^2(\widehat x-2)}{4}, ~\widehat\phi_3(\widehat x)=\frac{(\widehat x+1)^2(\widehat x-1)}{4}.$$ One can easily  verify that $\widehat{\phi}_0(-1)=1$,
$\widehat{\phi}_1'(-1)=1$, $\widehat{\phi}_2(1)=1$ and
$\widehat{\phi}_3'(1)=1$.

Legendre polynomials and Chebyshev polynomials are two most popular Jacobi polynomials.
Now we adopt Legendre polynomials $\{\widehat L_j\}_{j=0}^N$ or Chebyshev polynomials $\{\widehat T_j\}_{j=0}^N$  to give the bubble  basis functions on $I$.
One may set
$$\widehat\phi_{j}(\widehat x)=(2j-1)\widehat L_{j-4}(\widehat x) -2 (2j-3)\widehat L_{j-2}(\widehat x)+ (2j-5)\widehat L_{j}(\widehat x),~j=4,5,\cdots,N,$$
since $\widehat\phi_{j}(\widehat x)=\frac{(2j-1)(2j-3)(2j-5)}{4(j-2)(j-3)}\widehat{J}_{j}^{-2,-2}(\widehat x);$
another choice  is
$$\widehat\phi_{j}(\widehat x)=(j-1)\widehat T_{j-4}(\widehat x) -2 (j-2)\widehat T_{j-2}(\widehat x)+ (j-3)\widehat T_{j}(\widehat x),~j=4,5,\cdots,N.$$
\indent{\bf Remark 2.}~ One may set different polynomial degrees for each element. 
Figure 1 shows three element $\kappa_1-\kappa_3$ and the tensorial basis functions of $P_{4}(\kappa_1)$ on $\kappa_1$. If one wants to decrease the polynomial degrees to 3 on $\kappa_1$ and $\kappa_2$, one should delete the basis functions $\phi_4(x_1)\phi_4(x_2)$ on both $\kappa_1$ and $\kappa_2$ , $\phi_2(x_1)\phi_4(x_2)$ and $\phi_3(x_1)\phi_4(x_2)$ on the common edge of $\kappa_1$ and $\kappa_2$.

\begin{center}
\scalebox{1}{\includegraphics[width=1.1\textwidth]{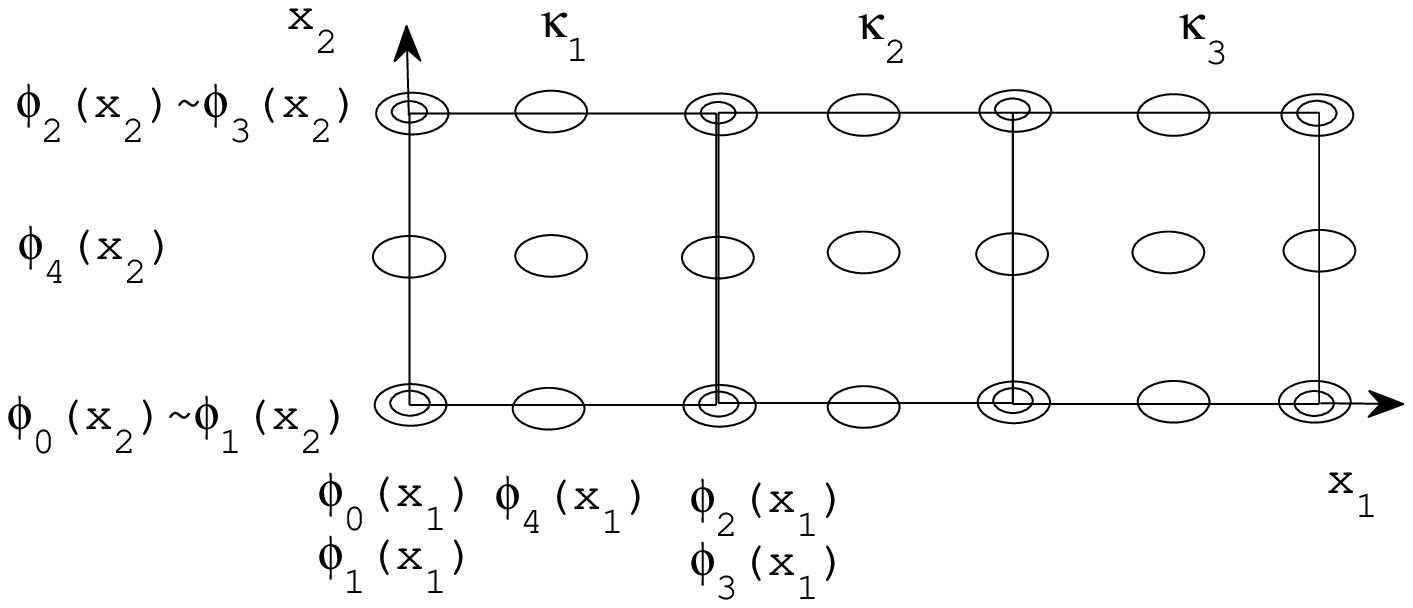}}\\\textrm{\small\bf Figure 1. The elements $\kappa_1-\kappa_3$ and the basis functions on   $\kappa_1$. }
\end{center}


\indent{\bf Remark 3.}~ The $H^m$-conforming spectral elements  can deal with the problem with mixed boundary condition on multi-dimensional domain
 due to adopting the nodal basis functions at the endpoint $\pm1$. Though we restrict our attention to the spectral   method on rectangular domain, it can be extended to   spectral   method on   non-rectangular   domains like the way as in \cite{guo}. Likewise the mesh adopted by the spectral element method  can be improved for  approximating
general domains better.

\section{$H^2$-conforming spectral elements for   transmission eigenvalues}
\label{sec:2}
In this section, we   aim  to apply the $H^2$-conforming spectral elements in the last section to the   transmission eigenvalue problem.

 %


 \indent Consider the Helmholtz transmission eigenvalue
problem: Find $k\in \mathcal{C}$, $\omega, \sigma\in L^{2}(D)$,
$\omega-\sigma\in H^{2}(D)$ such that
\begin{eqnarray}\label{s2.1}
&&\Delta \omega+k^{2}n(x)\omega=0,~~~in~ D,\\\label{s2.2}
 &&\Delta
\sigma+k^{2}\sigma=0,~~~in~ D,\\\label{s2.3}
 &&\omega-\sigma=0,~~~on~ \partial
D,\\\label{s2.4}
 &&\frac{\partial \omega}{\partial \nu}-\frac{\partial
\sigma}{\partial \nu}=0,~~~ on~\partial D,
\end{eqnarray}
where $D\subset \mathbb{R}^{d}~(d=2,3)$  is an open bounded simply
connected inhomogeneous medium,
 $\nu$ is the unit outward normal to $\partial D$.\\

 The
eigenvalue problem (\ref{s2.1})-(\ref{s2.4}) can be stated as
 the classical weak formulation below (see, e.g., \cite{rynne,cakoni2,colton2}): Find $k^{2}\in \mathcal{C}$, $k^{2}\not=0$, nontrivial $u\in
H_{0}^{2}(D)$ such that
\begin{eqnarray}\label{s2.5}
(\frac{1}{n(x)-1}(\Delta u+k^{2} u),\Delta v+\overline{k}^{2} n(x)v
)_{0}=0,~~~\forall v \in~H_{0}^{2}(D),
\end{eqnarray}
where $(\cdot,\cdot)_0$ is the inner product of $L^2(D)$.
  As usual, we  define $\lambda=k^{2}$ as   the transmission eigenvalue in this paper.
We suppose that the index of refraction $n\in
L^{\infty}(D)$ is a real valued function such that $n-1$ is strictly positive (strictly negative) almost everywhere in $D$.


\indent Define Hilbert space $\mathbf{H}=H_{0}^{2}(D)\times
L^{2}(D)$ 
and define
$\mathbf{H}^{s}(K)=H^{s}(K)\times H^{s-2}(K)$ with norm
$\|(u,w)\|_{s,K}=\|u\|_{s,K}+\|w\|_{s-2,K}$ for a given $K\subseteq D$.  We write $\mathbf H^1 := \mathbf H^1(D)$ for simplicity.\\

Although the problem (\ref{s2.5}) is a quadratic eigenvalue problem,  it can be linearized by introducing some variables.
 Using the linearized way   in \cite{yang1}, we introduce $w=\lambda u$,
then the problem (\ref{s2.5}) is equivalent to the following linear
weak formulation: Find $(\lambda, u, w)\in \mathcal{C}\times
H_{0}^{2}(D)\times L^{2}(D)$ such that
\begin{eqnarray}\label{2.6}
&&(\frac{1}{n-1}\Delta u, \Delta v)_{0}=\lambda
(\nabla(\frac{1}{n-1}u),
\nabla v)_{0}\nonumber\\
&&~~~~~~+\lambda(\nabla u, \nabla(\frac{n}{n-1}v))_{0}-\lambda
(\frac{n}{n-1}w, v )_{0},~~~\forall v
\in~H_{0}^{2}(D),\\ \label{2.7}
 &&(w,
z)_{0}=\lambda(u,z)_{0},~~~\forall z\in L^{2}(D).
\end{eqnarray}

\indent We introduce the following sesquilinear forms
\begin{eqnarray*}
&&A((u,w),(v,z))=(\frac{1}{n-1}\Delta u, \Delta v)_{0}+(w, z)_{0},\\
&&B((u,w),(v,z))\nonumber\\
&&~~~=(\nabla(\frac{1}{n-1}u), \nabla v)_{0}+(\nabla u,
\nabla(\frac{n}{n-1}v))_{0}- (\frac{n}{n-1}w, v
)_{0}+(u,z)_{0},
\end{eqnarray*}
then (\ref{s2.5})  can be rewritten as the following problem: Find
$\lambda\in \mathcal{C}$,  nontrivial $(u,w)\in \mathbf{H}$ such that
\begin{eqnarray}\label{2.8}
A((u,w),(v,z)) =\lambda B((u,w),(v,z)),~~~\forall (v,z)\in
\mathbf{H}.
\end{eqnarray}
Let norm $\|\cdot\|_A $ be induced by the inner product $A(\cdot,\cdot)$, then it is clear  $\| \cdot \|_A$ is equivalent to $\| \cdot \|_{2,D}$ in $\mathbf{H}$. \\




\indent One can easily verify that for any
given $(f,g)\in \mathbf{H}^{1}$, $B((f,g), (v, z))$ is a
continuous linear form on
$\mathbf{H}^1$:
\begin{eqnarray}\label{2.9}
  B((f,g), (v, z))\lesssim \|(f,g)\|_{1,D}\|(v,z)\|_{1,D},~\forall (v,z)\in \mathbf{H}^1.
\end{eqnarray}

\indent Consider the dual problem of (\ref{2.8}): Find
$\lambda^{*}\in \mathcal{C}$, nontrivial $(u^{*},w^{*})\in \mathbf{H}$
such that
\begin{eqnarray}\label{2.10}
A((v,z), (u^{*},w^{*})) =\overline{\lambda^{*}} B((v,z),
(u^{*},w^{*})),~~~\forall (v,z)\in \mathbf{H},
\end{eqnarray}
%
%
%
In order to discretize the space
$\mathbf{H}$, we need   finite element spaces to  discretize $ H_{0}^{2}(D)$ and
$L^{2}(D)$ respectively. Since $H^2_0(D)\subset L^2(D)$ here we can construct the  spectral element space $S_0^{N,h}:=S^{N,h}\cap H^2_0(D)$  such that $\mathbf{H}_{N,h}:=S_0^{N,h}\times S_0^{N,h}\subset H_{0}^{2}(D) \times L^{2}(D)$ . 

\indent The conforming spectral element approximation of (\ref{2.8})
is given by the following: Find $\lambda_{N,h}\in \mathcal{C}$, nontrivial
$(u_{N,h},w_{N,h})\in \mathbf{H}_{N,h}$ such that
\begin{eqnarray}\label{2.13}
A((u_{N,h},w_{N,h}),(v,z)) =\lambda_{N,h}
B((u_{N,h},w_{N,h}),(v,z)),~~~\forall (v,z)\in \mathbf{H}_{N,h}.
\end{eqnarray}

According to Theorem 2.4, we know the following error estimates hold for spectral element space.
For any $\psi\in H^2_0(D)\cap H^{2+r}(D)~(0\le r\le N-1)$  there holds
\begin{eqnarray*}
\inf\limits_{v\in S^{N,h}}\|\psi-v\|_s\lesssim (h/N)^{2+r-s}\|\psi\|_{2+r},~s=0,1,2.
\end{eqnarray*}

 \indent To give the error of eigenfunction $(u_{N,h},w_{N,h})$ in the norm $\|\cdot\|_{1,D}$ we need the following regularity assumption:\\
\indent {\bf R(D)}.~~ For any $\varrho\in H^{-1}(D)$, there
exists $\psi\in H^{2+r_{1}}(D)$ satisfying
\begin{eqnarray*}
&&\Delta(\frac{1}{n-1}\Delta \psi)=\varrho,~~~in~D,\\
&&\psi=\frac{\partial \psi}{\partial \nu}=0~~~ on~\partial D,
\end{eqnarray*}
and
\begin{eqnarray}\label{2.14}
\|\psi\|_{2+r_{1}}\leq C_{p} \|\varrho\|_{-1},
\end{eqnarray}
where $r_{1}\in (0,1]$, $C_{p}$ denotes the prior
constant dependent on the equation and $D$ but
independent of the right-hand side $\varrho$ of the equation.

It is easy to know that  (\ref{2.14})  is valid with $r_1 = 1$ when $n$ and $\partial D$ are appropriately smooth.
When $D \subset R^2$ is a convex polygon, from Theorem 2 in \cite{blum}
we can get $r_1 = 1$ if $n \in W^{2,p}(D)$. 
In this paper, let $\lambda$ be an eigenvalue of (\ref{2.8}) with the ascent $\alpha$. 
Let $M(\lambda)$ and $M(\lambda_{N,h})$ be the
space spanned by all generalized eigenfunctions corresponding respectively to the
eigenvalues $\lambda$ and $\lambda_{N,h}$. 
 As for the dual problem (\ref{2.10}), 
 the definitions
 of $M^*(\lambda^*)$ 
  are made similarly to $M(\lambda)$. 

In what follows, to characterize the approximation  of the finite element space $\mathbf{H}_{N,h}$ to 
$ M({\lambda})$ and $ M^*({\lambda^*})$, we introduce the following quantities
 \begin{eqnarray*}
 && {\delta}_{N,h}({{\lambda}})=\sup\limits_{(v,z)\in M({\lambda})\atop\|(v,z)\|_{2,D}=1}\inf\limits_{(v_h,z_h)\in \mathbf{H}_{N,h}}\|(v,z)-(v_{N,h},z_{N,h})\|_{2,D} , \\
 && {\delta}^*_{N,h}({{\lambda}}^*)=\sup\limits_{(v,z)\in M^*({\lambda}^*)\atop\|(v,z)\|_{2,D}=1}\inf\limits_{(v_{N,h},z_{N,h})\in \mathbf{H}_{N,h}}\|(v,z)-(v_{N,h},z_{N,h})\|_{2,D}.
 \end{eqnarray*}

Note that if $M(\lambda)\subset  \mathbf H^{t}(D)$ and $M^*(\lambda^*)\subset  \mathbf H^{t}(D)$ with
$t\le N+1$ then ${\delta}_{N,h}({{\lambda}})\lesssim
 (h/N)^{t-2} $ and ${\delta^*}_{N,h}({{\lambda^*}})\lesssim
 (h/N)^{t-2} $.

Using the spectral approximation theory \cite{babuska1,chatelin}, the literature \cite{yang1} established the following a priori error estimates for  the finite element approximation (\ref{2.13}).
According to \cite{yang1}, the following a priori error estimates are valid for the spectral elements,  as well as for   the spectral method as a special case.

\begin{theorem}
 Suppose $n \in W^{1,\infty}(D)$  and $\mathbf{R(D)}$ is valid. Let $\lambda_{N,h}$ be an eigenvalue of the  problem (\ref{2.13}) that converges to $\lambda$.
Let $(u_{N,h},w_{N,h})\in M(\lambda_{N,h})$
  and $\|(u_{N,h},w_{N,h})\|_A=1$, then there exists
$(u,w)\in M(\lambda)$ such that
\begin{eqnarray}\label{2.16}
&&\|(u_{N,h},w_{N,h})-(u,w)\|_{2,D}\lesssim
\delta_{N,h}(\lambda)^{1/\alpha} ,\\\label{2.17}
&&\|(u_{N,h},w_{N,h})-(u,w)\|_{1,D}\lesssim
(h/N)^{r_1/\alpha} \delta_{N,h}(\lambda)^{1/\alpha},\\\label{2.18}
&&|\lambda-\lambda_{N,h}|\lesssim
( \delta_{N,h}(\lambda)\delta_{N,h}^{*}(\lambda^{*}))^{1/\alpha}.
\end{eqnarray}
\end{theorem}

 \section{Numerical Experiment}
\indent In this section, we will report  some numerical
experiments for solving the transmission eigenvalue problem (\ref{2.8}) by the $H^2$ conforming spectral element method (SEM) on non-rectangular domain or by the spectral method (SM) on  rectangular domain.
Notice that the spectral scheme in \cite{an2} is based on the iterative method in \cite{sun1}, which is different from the one in this paper. An obvious feature of the method in \cite{an2} is using an estimated eigenvalue to initialize the iterative procedure.
We  consider the case when $D$ is the unit square or the L-shaped domain  in $d$-dimension ($d=2,3$)  and
the index of refraction $n=16,f_1(x),f_2(x)$ with $f_1(x)=8+x_{1}-x_{2},~f_2(x)=4+e^{x_1+x_2}$. We use uniform rectangular refinement to obtain some partitions  of $D$. Accordingly some numerical eigenvalues on the unit squares and the L-shaped are listed in Tables 1-5. We also depict the profiles   of some eigenfunctions on the L-shaped domains with $n=16$ (see Figure 2).\\
\indent We use Matlab 2012a to solve (\ref{2.13}) by the sparse solver  $eigs$ on a
Lenovo G480 PC with 4G memory. 
 For reading conveniently, we denote by $k_{j}=\sqrt{\lambda_{j}}$ and $k_{j,h}=\sqrt{\lambda_{j,h}}$ the $jth$   eigenvalue and the $jth$ numerical eigenvalue
obtained  on the space $\mathbf{H}_{h}$.


\indent  
Tables 1 and 3 show that the numerical eigenvalues obtained by SM  on the unit square in both two and three dimensions with different $n$ own superior accuracy; more precisely, they achieve about eight-digit accuracy with $N=15$.
Whereas the numerical eigenvalues obtained by SEM on the two L-shaped domains (see Tables 2, 4-5) do not have such accuracy.  This phenomenon is due to the eigenfunctions on the unit square are often smooth whereas those on the L-shaped domains have the singularities  at the L-corner point (see Figure 2).
\begin{figure}
\begin{tabular}{cc}
\includegraphics[width=0.4\textwidth]{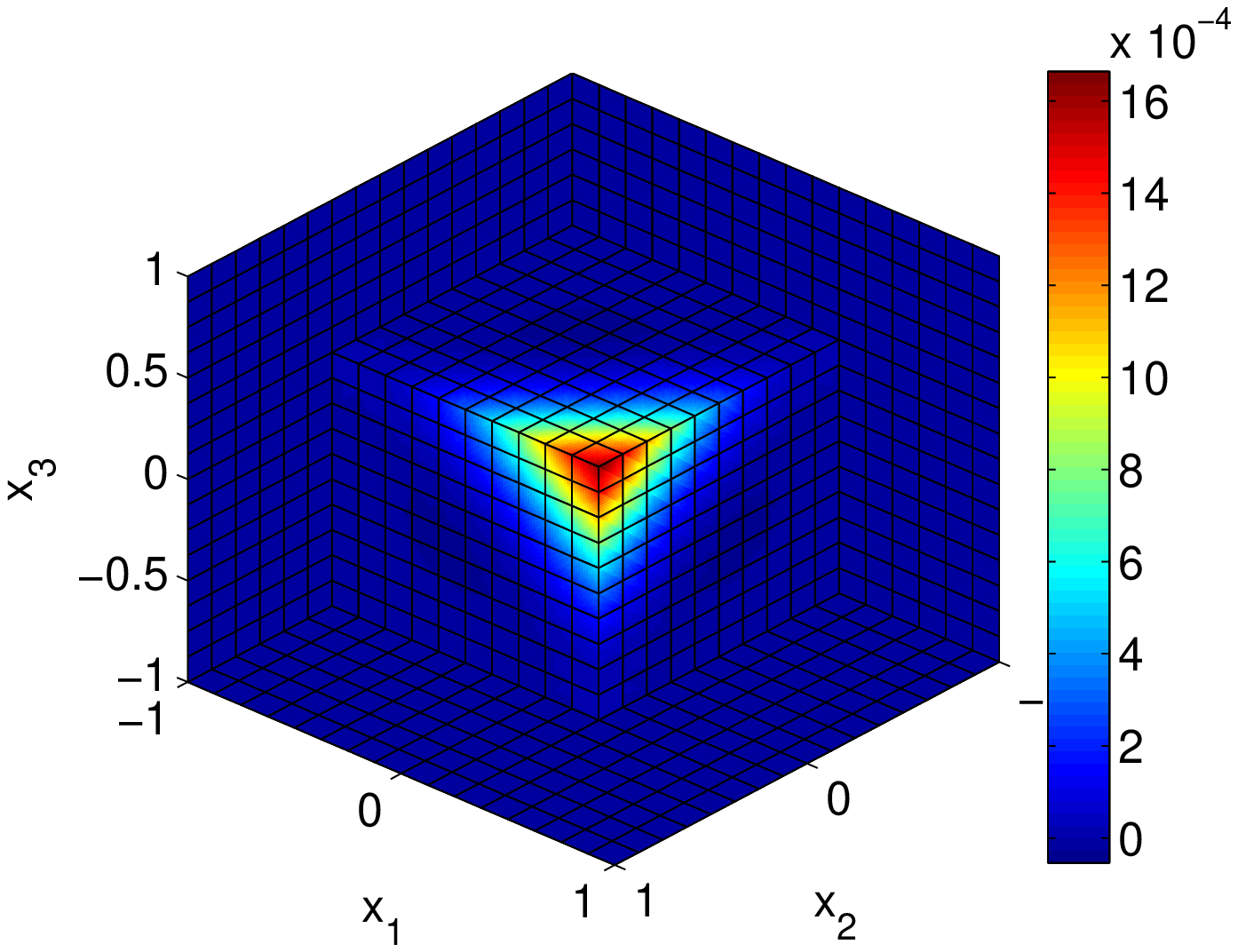}&
\includegraphics[width=0.4\textwidth]{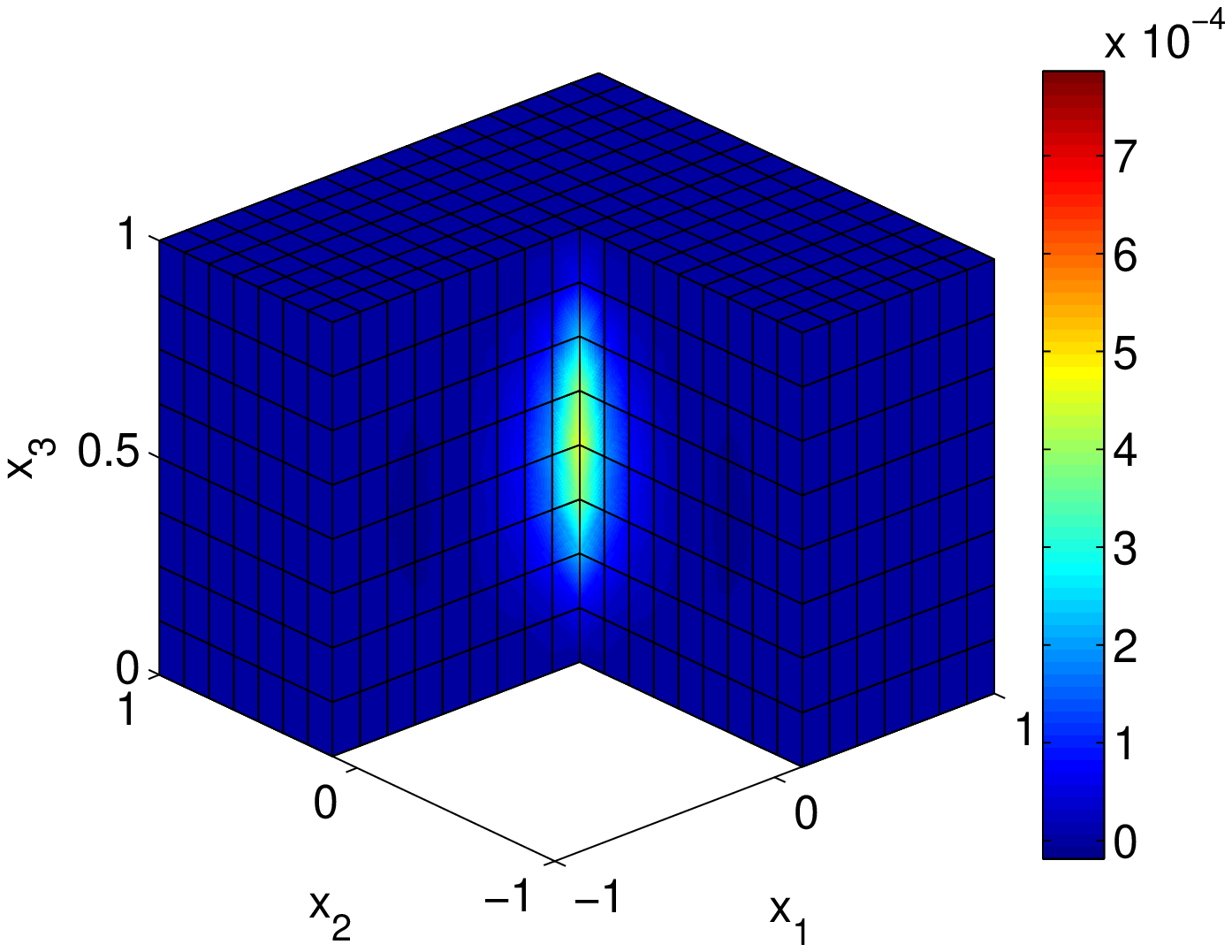}\\
\includegraphics[width=0.4\textwidth]{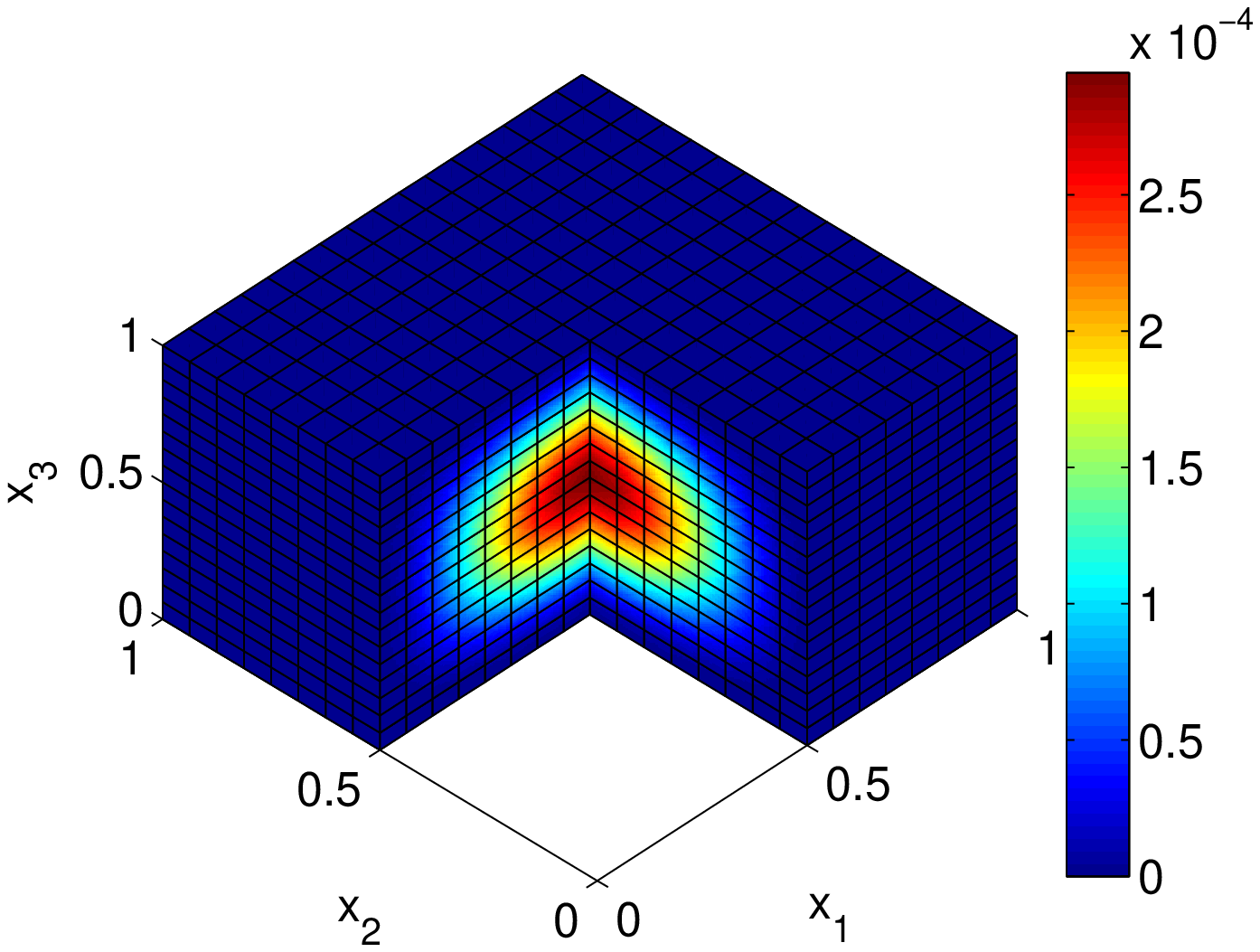}&
\includegraphics[width=0.4\textwidth]{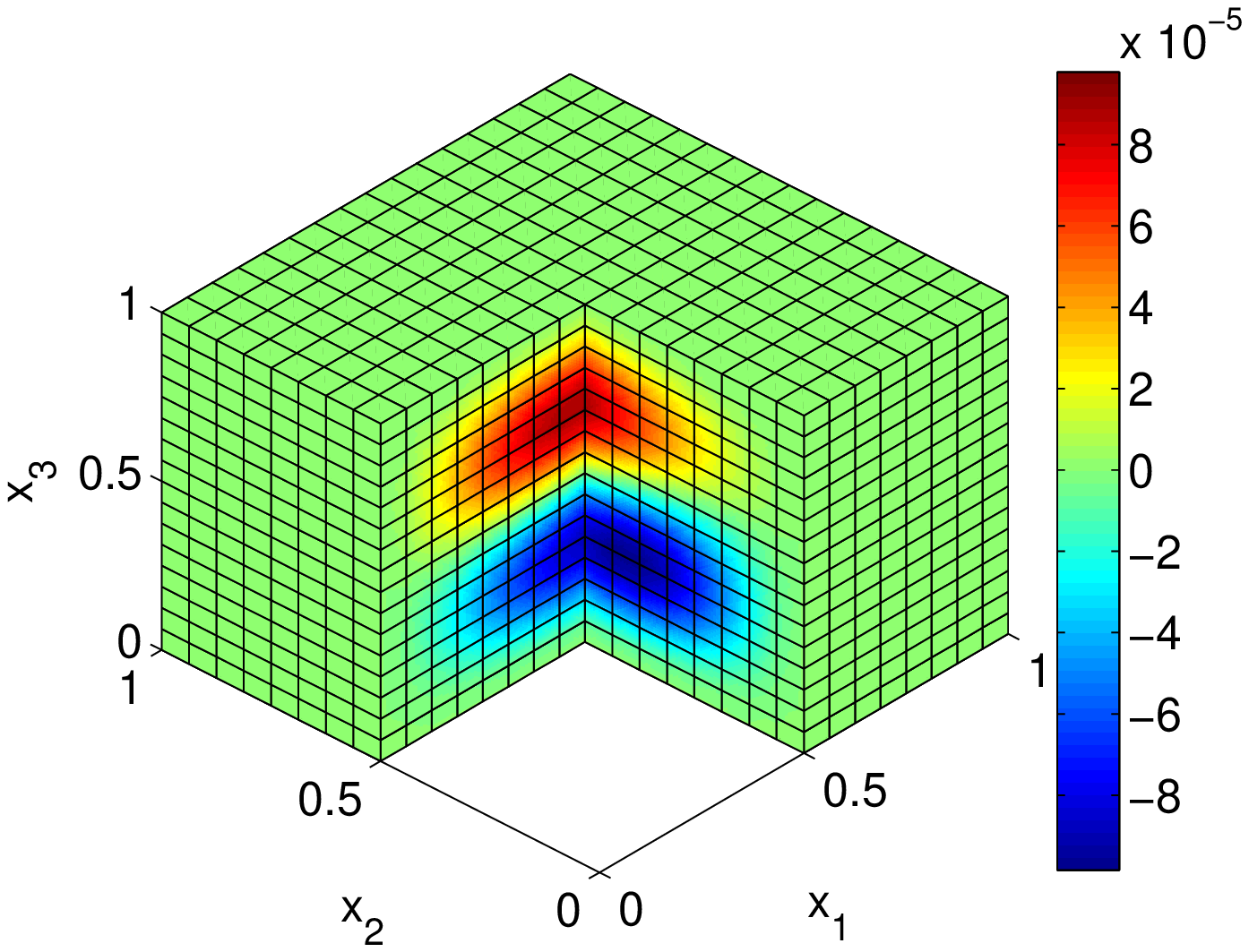}
\end{tabular}
\caption*{\textrm{\small\bf Figure 2. The 1st eigenfunction for $n=16$ on $(-1, 1)^3 \backslash (-1, 0)^3$   (left top), on  $((-1, 1)^2 \backslash (-1, 0]^2)\times(0, 1)$ (right top)  and on $(-1, 1)^3 $ (left bottom);
the 2nd eigenfunction for $n=16$ on $(-1, 1)^3$ (right bottom).  }}
\end{figure}

\begin{table}
\caption{Numerical eigenvalues obtained by SM on $(-0.5, 0.5)^2$.}
\begin{center} \footnotesize
\begin{tabular}{llllllll}\hline
$n$&$N$&$dof$&$k_{1,h}$&$k_{2,h}$,$k_{3,h}$&$k_{4,h}$&$k_{5,h}$\\\hline
 $16$&15&   288&      1.87959117836&   2.4442361007&   2.86643909864&   3.14011071773664\\
 $16$&20& 578&   1.87959117345&   2.4442360999&       2.86643911078&   3.14011071380238\\
 $16$&25& 968&     1.87959117325&   2.4442360993&     2.86643910989&   3.14011071380235\\
 $16$&30&  1458&    1.87959117313&    2.4442360992&   2.86643910981&   3.14011071380234\\
\hline\hline
$n$&$N$&$dof$&$k_{1,h}$&$k_{2,h}$&$k_{3,h}$&$k_{5,h},k_{6,h}$\\\hline
 $f_1$&15&   288& 2.8221893622&  3.5386966987&  3.5389915453&4.4965518722\\
 &&&&&&$\pm$0.8714816081i\\
 $f_1$&20& 578& 2.8221893415&  3.5386966965&  3.5389915430&  4.4965519547\\
 &&&&&&$\pm$0.8714817833i\\
 $f_1$&25& 968& 2.8221893411&  3.5386966953&  3.5389915418&  4.4965519545\\
 &&&&&&$\pm$0.8714817812i\\
 $f_1$&30&  1458&   2.8221893409&  3.5386966952&  3.5389915416&  4.4965519545\\
 &&&&&&$\pm$0.8714817805i\\
\hline\hline
$n$&$N$&$dof$&$k_{1,h}$,$k_{2,h}$&$k_{3,h}$&$k_{4,h}$&$k_{5,h}$\\\hline
 $f_2$&15&   288& 4.3184549937&  4.5885145655&  4.6472932515& 4.95760056967        \\
 &&&$\pm$0.6549618762i\\
 $f_2$&20& 578&   4.3184553572&  4.5885144838&  4.6472932378& 4.95759998805            \\
  &&&$\pm$0.6549618008i\\

 $f_2$&25& 968&   4.3184553557&  4.5885144805&  4.6472932351& 4.95759999028    \\
 &&&$\pm$0.6549617996i\\
 $f_2$&30&  1458& 4.3184553554&  4.5885144801&  4.6472932348& 4.95759999016       \\
 &&&$\pm$0.6549617990i\\
\hline
\end{tabular}
\end{center}
\end{table}

\begin{table}
\caption{Numerical eigenvalues obtained by SEM on $(-1, 1)^2\backslash([0, 1)\times(-1, 0]$).}
\begin{center} \footnotesize
\begin{tabular}{llllllll}\hline
$n$&$h$&$N$&$dof$&$k_{1,h}$&$k_{2,h}$&$k_{3,h}$&$k_{4,h}$\\\hline
 $16$&$\sqrt2$&15&   960&      1.47854&   1.569782&   1.705721&   1.78312049\\
 $16$&$\sqrt2$&20& 1870&    1.47742&   1.569746&   1.705408&   1.78311760\\
 $16$&$\sqrt2$&25&  3080&        1.47691&   1.569735&   1.705269&   1.78311674\\
 $16$&$\sqrt2$&30&   4590&      1.47665&   1.569730&   1.705195&   1.78311641\\
 $16$&$\frac{\sqrt2}{2}$&15& 4264&   1.47722&   1.569741&   1.705355&   1.78311725\\
 16&$\frac{\sqrt2}{4}$&15&  17928&    1.47663&   1.569730&   1.705189&   1.78311639\\
 $16$&$\frac{\sqrt2}{8}$&15&  73480&   1.47635&   1.569727&   1.705111&   1.78311614\\
\hline\hline
$n$&$h$&$N$&$dof$&$k_{1,h}$&$k_{2,h}$&$k_{3,h}$&$k_{5,h},k_{6,h}$\\\hline
 $f_1$&$\sqrt2$&15&   960&    2.30499&  2.395810&  2.64178134&  2.92613\\
 &&&&&&&$\pm$0.56694i\\
 $f_1$&$\sqrt2$&30&   4590&   2.30277&  2.395702&  2.64177929&  2.92466\\
 &&&&&&&$\pm$0.56511i\\
 $f_1$&$\frac{\sqrt2}{2}$&15& 4264&    2.30343&  2.395724&  2.64177970& 2.92510\\
 &&&&&&&$\pm$0.56566i\\
 $f_1$&$\frac{\sqrt2}{4}$&15&  17928&     2.30274&  2.395701&  2.64177927&  2.92464\\
 &&&&&&&$\pm$0.56509i\\
 $f_1$&$\frac{\sqrt2}{8}$&15&  73480&    2.30241&  2.395694&  2.64177916&  2.92442\\
  &&&&&&&$\pm$0.56482i\\
\hline
\end{tabular}
\end{center}
\end{table}

\begin{table}
\caption{Numerical eigenvalues obtained by SM on $(0,1)^3$.}
\begin{center} \footnotesize
\begin{tabular}{llllllll}\hline
$n$&$N$&$dof$&$k_{1,h}$&$k_{2,h},k_{3,h},k_{4,h}$&$k_{5,h},k_{6,h},k_{7,h}$&$k_{8,h},k_{9,h}$\\\hline
 16& 5& 8&        2.094055156&   2.664272514&3.0661457744& 3.406897998\\
 16& 10& 343&     2.067227464& 2.584867751&  2.9870636216& 3.246721378\\
 16& 15& 7304&   2.067227678&   2.584856761& 2.9870431376& 3.246569769\\
 16& 20& 4913&   2.067227671&   2.584856755& 2.9870431377& 3.246569737\\
\hline\hline
$n$&$N$&$dof$&$k_{1,h}$&$k_{2,h}$&$k_{3,h}$&$k_{4,h}$\\\hline
$f_1$&  5& 8&         3.098469114& 3.865559775&3.8745648277&3.877187844\\
$f_1$& 10& 343&       3.025670231&3.722083630&3.7247284048&3.724785357\\
$f_1$& 15& 7304&      3.025670590&3.722061785&3.7247087240&3.724765624\\
$f_1$& 20& 4913&    3.025670572&  3.722061777&3.7247087161&3.724765616\\
\hline
\end{tabular}
\end{center}
\end{table}

\begin{table}
\caption{Numerical eigenvalues obtained by SEM on $((-1, 1)^2 \backslash (-1, 0]^2)\times(0, 1)$.}
\begin{center} \footnotesize
\begin{tabular}{llllllll}\hline
$n$&$h$&$N$&$dof$&$k_{1,h}$&$k_{2,h}$&$k_{3,h}$&$k_{4,h}$\\\hline
 16&   $\sqrt3$& 4& 14& 1.85647&1.90219&1.96071&2.07418\\
 16& $\frac{\sqrt3}{2}$& 4& 512&  1.82025&   1.87691&   1.93480&   2.04285\\
 16& $\frac{\sqrt3}{4}$& 4& 6800&  1.80961&   1.87089&   1.92939&   2.03857\\
 16& $\sqrt3$& 7& 512&     1.81154&   1.87122&   1.92966&   2.03879\\
 16& $\sqrt3$& 8& 950& 1.80956&   1.87078&   1.92929&   2.03848\\
 16& $\sqrt3$& 9& 1584&    1.80841&   1.87069&   1.92925&   2.03845\\
 16& $\sqrt3$& 10& 2450&    1.80760&   1.87064&   1.92922&   2.03843\\

\hline\hline
$n$&$h$&$N$&$dof$&$k_{1,h}$&$k_{2,h}$&$k_{3,h}$&$k_{4,h}$\\\hline
$f_1$&  $\sqrt3$& 4& 14&  2.69988&  2.73831&  2.92444&  3.07209\\
$f_1$&  $\frac{\sqrt3}{2}$& 4& 512&      2.65535&   2.66530&   2.85574&   2.98575\\
$f_1$& $\frac{\sqrt3}{4}$& 4& 6800&        2.64344&   2.64820&   2.84249&   2.97109\\
$f_1$& $\sqrt3$& 7& 512&    2.64455&   2.65045&   2.84339&   2.97245\\
$f_1$& $\sqrt3$& 8& 950&    2.64325&   2.64810&   2.84230&   2.97092\\
$f_1$& $\sqrt3$& 9& 1584&      2.64269&   2.64702&   2.84203&   2.97042\\
$f_1$& $\sqrt3$& 10& 2450&  2.64215&   2.64641&   2.84186&   2.97005\\
\hline
\end{tabular}
\end{center}
\end{table}

\begin{table}
\caption{Numerical eigenvalues obtained by SEM on $(-1, 1)^3 \backslash (-1, 0)^3$.}
\begin{center} \footnotesize
\begin{tabular}{llllllll}\hline
$n$&$h$&$N$&$dof$&$k_{1,h}$&$k_{2,h},k_{3,h}$&$k_{4,h}$,$k_{5,h}$&$k_{6,h}$\\\hline
 16&   $\sqrt3$& 4& 74&       1.4993&   1.5552&       1.6676&1.6857\\
 16& $\frac{\sqrt3}{2}$& 4& 1568&       1.4443&   1.5199&    1.6443&   1.6515\\
 16& $\frac{\sqrt3}{4}$& 4& 17840&        1.4277&   1.5097&   1.6392&   1.6420\\
 16& $\frac{\sqrt3}{4}$& 5& 45808&    1.4225&   1.5069&   1.6388&   1.6395\\
 16& $\sqrt3$& 6& 774 &  1.4402&   1.5161&   1.6405&   1.6477  \\
 16& $\sqrt3$& 7&1568  & 1.4325&   1.5121& 1.6397&   1.6442   \\
 16& $\sqrt3$& 8&2770  &  1.4279&   1.5097& 1.6392&   1.6421  \\
 16& $\sqrt3$& 9& 4464&  1.4248&   1.5081  &  1.6389&  1.6406       \\
\hline\hline
$n$&$h$&$N$&$dof$&$k_{1,h}$&$k_{2,h}$&$k_{3,h}$&$k_{4,h}$\\\hline
 $f_1$&   $\sqrt3$& 4& 74&    2.2003&   2.2307&   2.2957&   2.3602\\
 $f_1$& $\frac{\sqrt3}{2}$& 4& 1568&     2.1191&   2.1675&   2.2282&   2.3227\\
 $f_1$& $\frac{\sqrt3}{4}$& 4& 17840&     2.0953&   2.1523&   2.2105&   2.3152\\
 $f_1$& $\frac{\sqrt3}{4}$& 5& 45808& 2.0880&   2.1480&   2.2066&   2.3143\\
 $f_1$& $\sqrt3$& 6& 774 & 2.1126&   2.1625&2.2197&   2.3178  \\
 $f_1$& $\sqrt3$& 7&1568  &   2.1020& 2.1561&   2.2141&2.3161   \\
 $f_1$& $\sqrt3$& 8&2770  & 2.0955&   2.1523&   2.2105&2.3151   \\
 $f_1$& $\sqrt3$& 9& 4464&     2.0913&   2.1498&   2.2083&   2.3146\\
\hline
\end{tabular}
\end{center}
\end{table}

\indent{\bf Acknowledgements.}~~
This work   is supported by the National Natural Science Foundation of China(Grant NO. 11561014).





\end{document}